\documentclass{amsart}
\usepackage{amssymb}
\usepackage{latexsym}
\usepackage{amscd}

\newtheorem{lem}{Lemma}
\newtheorem{conj}{Conjecture}
\newtheorem{prop}{Proposition}

\begin{document}
\renewcommand{\abstractname}{Abstract}
\begin{abstract}
In this note we prove a simple relation between the mean curvature
form, symplectic area, and the Maslov class of a Lagrangian immersion
in a K\"ahler-Einstein manifold. An immediate consequence is that in
K\"ahler-Einstein manifolds with positive scalar curvature, minimal
Lagrangian immersions are monotone.  
\end{abstract}
\title[A note on mean curvature ...]{A note on mean curvature, Maslov
  class and symplectic area  of Lagrangian immersions}     
\authors{Kai Cieliebak\footnote{kai@mathematik.uni-muenchen.de} and Edward Goldstein\footnote{egold@math.stanford.edu}}
\maketitle

\section{Introduction}
Let $(M,\omega)$ be a K\"ahler-Einstein manifold whose Ricci curvature
is a multiple of the metric by a real number $\lambda$. Hence the
K\"ahler form $\omega$ and the first Chern class $c_1(M)$ are related
by $c_1(M)=\lambda[\omega]$. Let $L$ be an immersed Lagrangian
submanifold of $M$ 
and let $\sigma_L$ be the mean curvature form of $L$ (which is a
closed 1-form on $L$). Let $F:\Sigma\to M$ be a smooth map from
a compact connected surface to $M$ whose boundary $\partial F$ is contained in
$L$. Let $\mu(F)$ be the Maslov class of $F$ and $\omega(F)$ its
symplectic area. The goal of this note is to prove the following
simple relation between these quantities:
\begin{equation}\label{eq:main}
   \mu(F) - 2\lambda\omega(F)=\frac{\sigma_L(\partial F)}{\pi}.
\end{equation}

This relation was given in~\cite{M} for $\mathbb{C}^n$ and in~\cite{Ar}
for Calabi-Yau manifolds. Dazord \cite{D} showed
that the differential of the mean curvature form is the Ricci form, so
in the K\"ahler-Einstein case $\sigma_L$ is closed. Y.G. Oh~\cite{Oh2}
investigated the symplectic area in the case that the mean curvature
form is exact.  

In the case $\lambda>0$, Lagrangian submanifolds for which the
left-hand side vanishes on all disks $F$ are called {\em monotone} in
the symplectic geometry literature, cf.~\cite{Oh1}. An immediate
consequence of (\ref{eq:main}) is that in K\"ahler-Einstein manifolds
with positive scalar curvature, minimal Lagrangian immersions are
monotone.  

In view of the condition $c_1(M)=\lambda[\omega]$, the left-hand side
of (\ref{eq:main}) depends only on the boundary of $F$. Thus if the map
$H_1(L;\mathbb{R})\to H_1(M;\mathbb{R})$ is trivial it defines a
cohomology class 
$\delta_L\in H^1(L;\mathbb{R})$ via $\delta_L(\gamma) := \mu(F) -
2\lambda\omega(F)$ for some 2-cycle $F$ with $\partial F=\gamma$. It follows
that in this case the cohomology class of the mean curvature form
$\sigma_L$ is invariant under symplectomorphisms of $M$. This
generalizes Oh's observation~\cite{Oh2} 
that the cohomology class is invariant under Hamiltonian
deformations. One consequence is the following:

Let $(M,\omega)$ be a K\"ahler manifold with $c_1(M)=\lambda
[\omega]\in H_2(M;\mathbb{R})$. Let $L$ be an immersed Lagrangian
submanifold of $M$ such that the map $H_1(L;\mathbb{R})\to
H_1(M;\mathbb{R})$ is trivial and $\delta_L\neq 0$. Suppose
there is a K\"ahler-Einstein metric $\omega_{KE}$ in the same
cohomology class as $\omega$ and $\phi:(M,\omega)
\rightarrow (M,\omega_{KE})$ is a symplectomorphism (e.g.~the one
provided by Moser's theorem). Then $\phi(L)$ is Lagrangian but {\em
not} minimal.  
 
Note that for $\lambda \neq 0$ most Lagrangian submanifolds $L$ with
nontrivial first Betti number such that  
$H_1(L)\to H_1(M)$ vanishes have $\delta_L\neq 0$: For any
such $L$, pick a normal vector field $v$ to $L$ such that $i_v
\omega$ is closed on $L$ and non-trivial cohomologically. Then small
time variations of $L$ through $v$ produce Lagrangian submanifolds
with nontrivial $\delta_L$.

\section{Notation}
We first recall the definition of the Maslov index that is suitable
for our purposes. Let $V$ be a Hermitian vector space of complex
dimension $n$. Let $\Lambda^{(n,0)}V$ be the (one-dimensional)
space of holomorphic $(n,0)$-forms on $V$ and set
$$
K^2(V):=\Lambda^{(n,0)}V\otimes\Lambda^{(n,0)}V.
$$
Let $L$ be a Lagrangian subspace of $V$. We can associate to $L$ an
element $\kappa(L)$ in $\Lambda^{(n,0)}V$ of unit length which
restricts to a real volume form on $L$. This element is unique up to
sign and therefore defines a unique element of unit length
$$
   \kappa^2(L):=\kappa(L)\otimes\kappa(L)\in K^2(V).
$$
Thus we get a map $\kappa^2$ from the Grassmanian $Gr_{\rm Lag}(V)$ of
Lagrangian planes to the unit circle in $K^2(V)$. This map induces
a homomorphism $\kappa^2_{\ast}$ of fundamental groups 
\[
   \kappa^2_{\ast}:\pi_1(Gr_{\rm Lag}(V)) \rightarrow \mathbb{Z}.
\] 
To understand the map $\kappa^2_{\ast}$, let $L$ be a Lagrangian subspace
and let $v_1,\ldots,v_n$ be an orthonormal basis for $L$. For $0 \leq
t \leq 1$ consider the subspace
\[
   L_t=span\{v_1,\ldots,v_{n-1},e^{\pi i t}v_n\}.
\] 
This loop $\{L_t\}$ is the standard generator of $\pi_1 (Gr_{\rm
  Lag}(V))$. The induced elements in $\Lambda^{(n,0)}V$ are related by
$\kappa(L_t)=\pm e^{-\pi i t}\kappa(L)$, so $\kappa^2(L_t)=e^{-2\pi i
  t}\kappa^2(L)$ and $\kappa^2_{\ast}(\{L_t\})=-1$. Thus we see that the 
homomorphism $\kappa^2_{\ast}$ is related to the Maslov index $\mu$ (as
defined, e.g., in~\cite{ALP}) by
$$
   \kappa^2_\ast = -\mu : \pi_1(Gr_{\rm Lag}(V)) \rightarrow \mathbb{Z}.
$$  
Now let $(M,\omega)$ be a symplectic manifold of dimension $2n$. Pick 
a compatible almost complex structure $J$ on $M$  and let $K(M)$ 
be the canonical bundle of $M$, i.e., $K(M):=\Lambda^{(n,0)}T^{\ast}M$ 
is the bundle of $(n,0)$-forms on $M$. Note that
$c_1(K(M))=-c_1(M)$. Let $K^2(M):=K(M)\otimes K(M)$ be the square of
the canonical bundle. \\  
Let $L$ be an immersed Lagrangian submanifold of $M$. For any point $l
\in L$ there is an element of unit length $\kappa(l)$ of $K(M)$ over
$l$, unique up to sign, which restricts to a real volume form on the
tangent space $T_lL$. The squares of these elements give rise to a
section of unit length
$$
\kappa^2_L:L \rightarrow K^2(M).
$$
Now let $F:\Sigma\to M$ be a smooth map with boundary $\partial F$ on
$L$. The {\em symplectic area} of $F$ is 
\[
   \omega(F)=\int_{\Sigma}F^*\omega.
\] 
This defines a map from the relative second homology group to
$\mathbb{R}$, 
$$
   [\omega]: H_2(M,L;\mathbb{Z}) \rightarrow \mathbb{R}.
$$
To define the Maslov class $\mu(F)$, choose a unitary frame for
the tangent bundle $TM$ along $F$. Consider the dual frame and
wedge all its elements. Thus we get a unit 
length section $\kappa_F$ of $K(M)$ over $F$. Now on the boundary
$\partial F=F(\partial\Sigma)$ we also have the section $\kappa^2_L$
defined above. We can uniquely write
$$
   \kappa_L^2 = e^{i\theta} \kappa_F^2
$$
for a function $e^{i\theta}:\partial\Sigma\to S^1$ to the unit
circle. The Maslov class $\mu(F)$ is minus its winding number,
$$
   \mu(F) := \frac{-1}{2\pi} \int_{\partial F}d\theta.
$$
This defines a map
$$
   \mu: H_2(M,L;\mathbb{Z}) \rightarrow \mathbb{Z}.
$$
In view of the discussion above, this definition agrees with the usual
definition of the Maslov class, cf.~\cite{ALP}. 

Now suppose that $c_1(M)=\lambda[\omega]\in H_2(M;\mathbb{R})$. 
It is well-known that if $\partial F$ is trivial in
$H_1(L;\mathbb{R})$, then $F$ represents an element in $[F]\in
H_2(M;\mathbb{R})$ and  
$$
   \mu(F) = 2c_1(M)([F]) = 2\lambda\omega(F).
$$
So in this case $\mu(F)-2\lambda\omega(F)$ depends only on the
boundary $\partial F\in H_1(L;\mathbb{R})$. If, moreover, the map
$H_1(L;\mathbb{R})\to H_1(M;\mathbb{R})$ is trivial, this expression defines a
cohomology class 
$\delta_L\in H^1(L;\mathbb{R})$ via 
$$
   \delta_L(\gamma) := \mu(F) - 2\lambda\omega(F)
$$ 
for some 2-cycle $F$ with $\partial F=\gamma$.

\section{proof}
Now assume that $(M,\omega)$ is K\"ahler-Einstein, i.e., $M$ carries a
K\"ahler metric whose Ricci curvature is a multiple of the metric by a
constant $\lambda\in\mathbb{R}$. This is equivalent to saying that the
curvature form of the canonical bundle $K(M)$ equals
$-\frac{2\pi}{i}\lambda\omega$. We denote the connections on $K(M)$
and $K^2(M)$ (induced by the Levi-Civita connection) by $\nabla$. \\ 
Let $L$ be an immersed Lagrangian submanifold of $M$ and let
$\kappa^2_L$ be the canonical section of $K^2(M)$ over $L$ as
above. The section $\kappa^2_L$ defines a connection 1-form $\xi_L$
for $K^2(M)$ over $L$ by the condition $\nabla\kappa^2_L = \xi_L
\otimes \kappa^2_L$.  Since $\kappa^2_L$ has constant length $1$,
$\xi_L$ is an imaginary valued 1-form on $L$. From the Einstein
condition and the fact that $L$ is Lagrangian we get $d(i\xi_L) =
-4\pi\lambda\omega|_L=0$, so the form $i\xi_L$ is closed. \\ 
Let $H$ be the trace of the second fundamental form of $L$ (the mean
curvature vector field of $L$). Thus $H$ is a section of the 
normal bundle to $L$ in $M$ and we have a corresponding 1-form
$\sigma_L:= i_H\omega$ on $L$. The following fact goes back to~
\cite{Oh2} (see also~\cite{Gold1} for a proof): 
$$
   \sigma_L= i \xi_L/2.
$$
(Here the factor $1/2$ is due to the fact that $\xi_L$ is a connection
1-form for $K^2(M)$ rather than $K(M)$.) Thus $\sigma_L$ is a closed
1-form on $L$, called the {\em mean curvature form} on $L$.   

Having explained all the terms in formula (\ref{eq:main}), we now
turn to its proof. Let $F:\Sigma\to M$ be a smooth map from a compact
surface with boundary on $L$. Define the section $\kappa_F$ of $K(M)$
over $F$ as above, using a unitary trivialization of $TM$ over $F$. 
Let $\xi_F$ be the connection 1-form along $F$ defined by 
$\nabla\kappa^2_F = \xi_F \otimes \kappa^2_F$. The Einstein
condition tells us that $d(i\xi_F) = -4\pi\lambda F^*\omega$. Thus by
Stokes' theorem, 
\[
   2\lambda\omega(F)=\int_{\partial F}\frac{-i\xi_F}{2\pi}.
\] 
Recall that along $\partial F$ we have
$\kappa^2_L=e^{i\theta}\kappa^2_F$ for a function
$e^{i\theta}:\partial\Sigma\to S^1$, and the Maslov class is given by 
$$
   \mu(F) = -\frac{-1}{2\pi} \int_{\partial F}d\theta.
$$
The connection 1-forms $\xi_F$ and $\xi_L$ are related by
\[
   \xi_L=\xi_F+i\,d\theta. 
\] 
Thus
\[
   \frac{\sigma_L(\partial F)}{\pi} = \int_{\partial
   F}\frac{i\xi_L}{2\pi} = \int_{\partial F}\frac{i\xi_F}{2\pi} -
   \int_{\partial F}\frac{d\theta}{2\pi} = \mu(F) - 2\lambda\omega(F).  
\]


\end{document}